\documentclass[letterpaper, 10 pt, conference]{ieeeconf}%
\usepackage{makeidx}
\usepackage{graphicx}
\usepackage{multicol}
\usepackage[bottom]{footmisc}
\usepackage{latexsym}
\usepackage{subfigure}
\usepackage{cite,color,comment,xspace}
\usepackage{amsfonts}
\usepackage{amsmath}
\usepackage{amssymb}
\usepackage{algorithm}
\usepackage{algorithmic}
\usepackage{url}
\usepackage{times}
\usepackage{mathptmx}
\usepackage{enumerate}
\usepackage{epstopdf}
\usepackage{subfigure}
\usepackage{epsfig}
\usepackage{authblk}
\usepackage{mathrsfs}
\usepackage[skip=2pt,font=scriptsize]{caption}%
\providecommand{\U}[1]{\protect\rule{.1in}{.1in}}
\renewcommand{\baselinestretch}{.985}
\IEEEoverridecommandlockouts

\makeatletter

\@addtoreset{corollary}{section}

\@addtoreset{lemma}{section}

\@addtoreset{definition}{section}
\makeatother
\allowdisplaybreaks
\begin{document}

\title{{\LARGE \textbf{Optimal Dynamic Formation Control of Multi-Agent Systems in
Environments with Obstacles}}}
\author{Xinmiao Sun and Christos G. Cassandras\thanks{{\footnotesize The authors' work
is supported in part by NSF under grants CNS-1239021 and IIP-1430145, by AFOSR
under grant FA9550-12-1-0113, and by ONR under grant N00014-09-1-1051.}%
}\thanks{{\footnotesize Division of Systems Engineering and Center for
Information and Systems Engineering, Boston University; e-mail: \{xmsun,cgc\}
@bu.edu}}}
\maketitle

\begin{abstract}
We address the optimal dynamic formation problem in mobile leader-follower
networks where an optimal formation is generated to maximize a given objective
function while continuously preserving connectivity. We show that in a convex
mission space, the connectivity constraints can be satisfied by any feasible
solution to a mixed integer nonlinear optimization problem. When the optimal
formation objective is to maximize coverage in a mission space cluttered with
obstacles, we separate the process into intervals with no obstacles detected
and intervals where one or more obstacles are detected. In the latter case, we
propose a minimum-effort reconfiguration approach for the formation which
still optimizes the objective function while avoiding the obstacles and
ensuring connectivity. We include simulation results illustrating this dynamic
formation process.

\end{abstract}

\section{Introduction}

The multi-agent system framework consists of a team of autonomous agents
cooperating to carry out complex tasks within a given environment that is
potentially highly dynamic, hazardous, and even adversarial. The overall
objective of the system may be time-varying and combines exploration, data
collection, and tracking to define a \textquotedblleft
mission\textquotedblright. Related problems are often referred to as
multi-agent coordination
\cite{Cao2013TIICoordination,boutilier1999sequential,beard2001coordination} or
cooperative control \cite{Shamma:2008,Choi20092802,cgc2005}. In many cases,
mobile agents are required to establish and maintain a certain spatial
configuration, leading to a variety of \emph{formation control} problems.
These problems are generally approached in two ways: in the leader-follower
setting, an agent is designated as a team leader moving on some given
trajectory with the remaining agents tracking this trajectory while
maintaining the formation; in the leaderless setting the formation must be
maintained without any such benefit. Examples of formation control problems
may be found in \cite{cao2011maintaining},\cite{Kwang2014}%
,\cite{Desai1999,Yamaguchi1994,jiannan2013,ji2007distributed} and references
therein. In robotics, this is a well-studied problem; for instance in
\cite{Yamaguchi1994}, a desired shape for a networked\ strongly connected
group of robots is achieved by designing a quadratic spread potential field on
a relative distance space. In \cite{Desai1999}, a leader and several followers
move in an area with obstacles which necessitate the transition from an
initial formation shape to a desired new shape; however, the actual choice of
formations for a particular mission is not addressed in \cite{Desai1999}, an
issue which is central to our approach in this paper. In
\cite{ji2007distributed} the authors consider the problem of preserving
connectivity when the nodes have limited sensing and communication ranges;
this is accomplished through a control law based on the gradient of an
edge-tension function. More recently, in \cite{jiannan2013}, the goal is to
integrate formation control with trajectory tracking and obstacle avoidance
using an optimal control framework.

In this paper, we take a different viewpoint of formations. Since agent teams
are typically assigned a mission, there is an objective (or cost) function
associated with the team's operation which depends on the spatial
configuration (formation) of the team. Therefore, we view a formation as the
result of an optimization problem which the agent team solves in either
centralized or distributed manner. We adopt a leader-follower approach,
whereby the leader moves according to a trajectory that only he/she controls.
During the mission, the formation is preserved or must adapt if the mission
(hence the objective function) changes or if the composition of the team is
altered (by additions or subtractions of agents) or if the team encounters
obstacles which must be avoided. In the latter case in particular, we expect
that the team adapts to a new formation which still seeks to optimize an
objective function so as to continue the team's mission by attaining the best
possible performance. The problem is complicated by the fact that such
adaptation must take place in real time. Thus, if the optimization problem
determining the optimal formation is computationally demanding, we must seek a
fast and efficient control approach which yields possibly suboptimal
formations, but guarantees that the initial connectivity attained is
preserved. Obviously, once obstacles are cleared, the team is expected to
return to its nominal optimal formation.

Although the optimal dynamic formation control framework proposed here is not
limited by the choice of tasks assigned to the team, we will focus on the
coverage control problem because it is well studied and amenable to efficient
distributed optimization methods
\cite{SM2011,CM2004,cgc2005,caicedo2008coverage,caicedo2008performing,breitenmoser2010voronoi,Minyi2011,Gusrialdi2011}%
, while also presenting the challenge of being generally non-convex and
sensitive to the agent locations during the execution of a mission. The local
optimality issue, which depends on the choice of objective function, is
addressed in \cite{Sun2014,schwager2008,gusrialdi2013improved}, while the
problem of connectivity preservation in view of limited communication ranges
is considered in \cite{ji2007distributed,Minyi2011}.

The contribution of this paper is to formulate an optimization problem which
jointly seeks to position agents in a two-dimensional mission space so as to
optimize a given objective function while at the same time ensuring that the
leader and remaining agents maintain a connected graph dictated by minimum
distances between agents, thus resulting in an \emph{optimal} formation. The
minimum distances may capture limited communication ranges as well as any
other constraint imposed on the team. We show that the solution to this
problem guarantees this connectivity. The formation becomes \emph{dynamic} as
soon as the leader starts moving along a given trajectory which may either be
known to all agents in advance or determined only by the leader. Thus, it is
the team's responsibility to maintain an optimal formation. We show that this
is relatively simple as long as no obstacles are encountered. When one or more
obstacles are encountered (i.e., they come within the sensing range of one or
more agents), then we propose a scheme for adapting with minimal effort  to a
new formation which maintains connectivity while still seeking to optimize the
original team objective.

The paper is organized as follows. In Sec. II, we formulate a general optimal
formation control problem. In Sec. III, we focus on a convex feasible space
and derive a mixed integer nonlinear problem whose solution is shown to ensure
connectivity while maintaining an optimal formation. In Sec. IV, we propose a
scheme to solve the optimal formation problem in a mission space with
obstacles. We propose an algorithm to first obtain a connected formation and
then optimize it while maintaining connectivity. Simulation results are
included in Sec. V.

\section{Optimal Formation Problem Formulation}

Consider a set of $N+1$ agents with a leader labeled $0$ and $N$ followers
labeled 1 through $N$ in a mission space $\Omega\in\mathbb{R}^{2}$. Agent $i$
is located at $s_{i}(t)\in\mathbb{R}^{2}$ and let $\mathbf{s}(t)=(s_{0}%
(t),...,s_{N}(t))$ be the full agent location vector at $t$. The leader
follows a predefined trajectory $s_{0}(t)$ over $t\in\lbrack0,T]$ which is
generally not known in advance by the remaining agents. We model the agent
team as a undirected graph
$\mathscr{G(\mathbf{s})}=(\mathscr{N},\mathscr{E},\mathbf{s})$, where
$\mathscr{N}=\{0,1,...,N\}$ is the set of agent indices and let
$\mathscr{N}_{F}=\{1,\ldots,N\}\subset\mathscr{N}$ be the set of follower
indices. In this model, the set of edges $\mathscr{E}=\{(i,j):i,j\in
\mathscr{N}\}$ contains all possible agent pairs for which constraints may be imposed.

In performing a mission, let $H(\mathbf{s}(t))$ be an objective function
dependent on the agent locations $\mathbf{s}(t)$. If the locations are
unconstrained, the problem is posed as $\max_{\mathbf{s}(t)\in\Omega
}H(\mathbf{s}(t))$ subject to dynamics that may characterize the motion of
each agent. If $t$ is fixed, then this is a nonlinear parametric optimization
problem over the mission space $\Omega$ \cite{Minyi2011}. If, on the other
hand, agents are required to also satisfy some constraints relative to each
other's position, then a \emph{formation} is defined as a graph that satisfies
these constraints. We then introduce a Boolean variable $c(s_{i},s_{j})$ to
indicate whether two agents satisfy these constraints:
\begin{equation}
{\label{connnectedCon}}c(s_{i},s_{j})=\left\{
\begin{array}
[c]{cl}%
1 & \text{all constraints are satisfied}\\
0 & \text{otherwise}%
\end{array}
\right.
\end{equation}
and if $c(s_{i},s_{j})=1$ we say that agents $i$ and $j$ are \emph{connected}.
A loop-free path from $i$ to the leader, $\pi_{i}=\{0,\ldots,a,b,\ldots,i\}$,
is defined as an ordered set where neighboring agents are connected such that
$c(s_{a},s_{b})=1$. Let $\Pi_{i}$ be the set of all possible paths connected
to the leader. The graph $\mathscr{G(\textbf{s})}$ is connected if $\Pi
_{i}\neq\emptyset$ for all $i\in\mathscr{N}_{F}$. We can now formulate an
optimal formation problem with connectivity preservation as follows, for any
fixed $t\in\lbrack0,T]$:
\begin{equation}%
\begin{split}
&  \max_{\mathbf{s}(t)\in\Omega}\text{ }H(\mathbf{s}(t))\\
\text{s.t.}\quad s_{i}(t)  &  \in F\subseteq\Omega,\text{ \ }i\in
\mathscr{N}_{F}\\
s_{0}(t)  &  \text{ is given}\\
\mathscr{G}(\mathbf{s}(t))  &  \text{ is connected}%
\end{split}
\label{obj1}%
\end{equation}
For the sake of generality, we impose the constraint $s_{i}(t)\in
F\subseteq\Omega$ for all follower agents to capture the possibility that a
formation is constrained. The \emph{feasible space} $F$ can be convex (e.g.,
followers may be required to be located on one side of the leader relative to
a line in $\Omega$ that goes through $s_{0}(t)$) or non-convex (e.g.,
followers may be forbidden to enter polygonal obstacles and $F$ is the set
$\Omega$ excluding all interior points of the obstacles). The solution to this
problem is an \emph{optimal formation} at time $t$ and is denoted by
$\mathscr{G}_{F}(\mathbf{s}(t))$. Given a time interval $[t_{1},t_{2}]$, the
formation is \emph{maintained} in $[t_{1},t_{2}]$ if $s_{i}(t)-s_{i}%
(t_{1})=s_{0}(t)-s_{0}(t_{1})$ holds for all $t\in\lbrack t_{1},t_{2}]$,
$i\in\mathscr{N}_{F}$; otherwise, it is a new formation. Figure
\ref{fig:missionSpace} shows an example of optimal dynamic formation control
in a mission space with obstacles. \begin{figure}[ptb]
\centering
\includegraphics[
height=1.5in,
width=3in]
{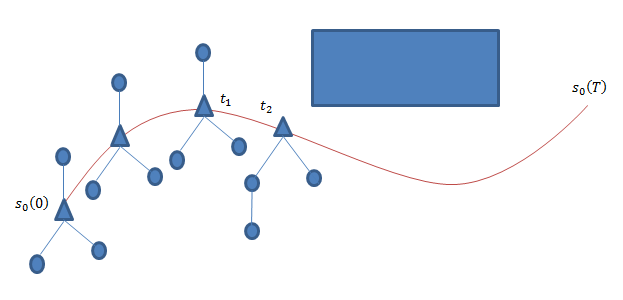} \caption{A mission space example where the triangle is
the leader and the red line is a predefined trajectory in $[0,T]$. The circles
are followers and the rectangle is an obstacle. The formation is maintained in
$[0,t_{1}]$, but at $t_{2}$ a new formation is needed.}%
\label{figce}%
\end{figure}Clearly, this is a challenging problem. To begin with, the last
constraint in (\ref{obj1}) is imprecise and may be different in a convex or
non-convex feasible space. In addition, the computational complexity of
obtaining a solution may be manageable in determining an initial formation but
becomes infeasible if a new formation $\mathscr{G}_{F}(\mathbf{s}(t))$ is
required during the real-time execution of a mission. In the following two
sections, we first propose an approach to solve this problem in a convex
feasible space and then use this solution to enable the maintenance of a
formation in a non-convex case.

\section{Optimal Dynamic Formation Control in a Convex Feasible Space}

In a convex feasible space, the simplest connection constraints are of the
form $d_{ij}(t)\equiv\Vert s_{i}(t)-s_{j}(t)\Vert\leq C_{ij}$ for some pair
$(i,j)$, $i,j\in\{0,1,...,N\}$, where $C_{ij}>0$ is a given scalar. This may
be the minimum distance needed to establish communication or $d_{ij}$ may be
used to enforce a specific desired shape in the formation. Techniques based on
the graph Laplacian \cite{merris1994laplacian} are often used to solve this
kind of problem, e.g., \cite{R2004}. However, our goal is to determine a
formation which solves the optimization problem in (\ref{obj1}) for a given
$H(\mathbf{s}(t))$. Thus, we describe next an approach to transform the last
constraint in (\ref{obj1}) into a mixed integer nonlinear optimization problem
by introducing a set of \emph{flow variables} over $\mathscr{G(\textbf{s})}$.
The leader $0$ is assumed to be a source node which sends $N$ units of flow
through the graph $\mathscr{G(\textbf{s})}$ to all other agents. Let
$\rho_{ij}\in\mathbb{Z}^{+},$ $i\in\mathscr{N},$ $j\in\mathscr{N}_{F}$ be an
integer flow amount through link $(i,j)$. Note that, in general, $\rho
_{ij}\neq\rho_{ji}$ and that either $\rho_{ij}>0$ or $\rho_{ji}>0$ implies
that $c(s_{i},s_{j})=1$. We can then define a flow vector $\mathbf{\rho}%
=(\rho_{01},\rho_{11},\rho_{N1},\ldots,\rho_{0N},\ldots,\rho_{NN})$. Observe
that $\rho_{i0},i\in\mathscr{N}$ is not a flow variable in $\mathbf{\rho}$
since the leader is not allowed to receive any flows from the followers. For
each follower $j$, we define an auxiliary variable $N_{j}$ to be the net flow
at node $j$:
\begin{equation}
N_{j}=\sum_{i\in\mathscr{N}}\rho_{ij}-\sum_{i\in\mathscr{N}_{F}}\rho_{ji}%
\end{equation}
Using this notation, we introduce next a number of linear constraints that
represent a connected graph. First, the leader provides $N$ units of flow:
\begin{equation}
\sum_{i\in\mathscr{N}_{F}}\rho_{0i}=N\label{leader_total}%
\end{equation}
Next, each follower $j$ must receive a net flow $N_{j}=1$ in order to ensure
that there is one path from the leader to $j$:
\begin{equation}
N_{j}=\sum_{i\in\mathscr{N}}\rho_{ij}-\sum_{i\in\mathscr{N}_{F}}\rho
_{ji}=1,\text{ \ }j\in\mathscr{N}_{F}\label{follower_conservation}%
\end{equation}
To prohibit self loops we require that
\begin{equation}
\rho_{ii}=0,\text{ \ }i\in\mathscr{N}\label{self_loop}%
\end{equation}
Finally, the maximal flow capacity is upper bounded by the source amount $N$:
\begin{equation}
\rho_{ij}\leq N,\text{ \ }i\in\mathscr{N},\text{ \ }j\in\mathscr{N}_{F}%
\label{flow_capacity}%
\end{equation}
Observe that (\ref{leader_total}) and (\ref{follower_conservation}) are
linearly dependent since $\sum_{j}N_{j}=N$. Thus, the constraint
(\ref{leader_total}) is redundant and may be omitted.

\textbf{Theorem 1} If there exists a flow vector $\mathbf{\rho}$ such that
constraints (\ref{follower_conservation})-(\ref{flow_capacity}) hold, then
there exists a connected graph $\mathscr{G(\mathbf{s})}$. Moreover, the number
of possible graphs is finite.

\textbf{Proof:} We use a contradiction argument. Assume that at least one
follower agent is not connected to the leader while satisfying
(\ref{follower_conservation})-(\ref{flow_capacity}). We can separate the
follower agents into two sets: $N_{1}=\{k:\Pi_{k}\neq\emptyset\}$ and
$N_{2}=\{j:\Pi_{j}=\emptyset\}$. Then, $\rho_{kj}=0$ must be true for all
$k\in N_{1}$ and $j\in N_{2}$. This is because if $\rho_{kj}>0$, then there
exists a path $\pi_{j}=\{\pi_{k},j\}$ where $\pi_{k}\in\Pi_{k}$, which
contradicts the fact that $j\in N_{2}$. In addition, obviously $\rho_{0j}=0$
for $j\in N_{2}$. Summing the left-hand-sides of all constraints
(\ref{follower_conservation}) such that $j\in N_{2}$, we obtain%
\begin{equation}%
\begin{split}
&  \sum_{j\in N_{2}}N_{j}=\sum_{j\in N_{2}}\left(  \sum_{k\in\mathscr{N}}%
\rho_{kj}-\sum_{k\in\mathscr{N}_{F}}\rho_{jk}\right)  \\
= &  \sum_{j\in N_{2}}\left[  \sum_{k\in N_{1}}\rho_{kj}+\sum_{k\in N_{2}}%
\rho_{kj}+\rho_{0j}-\left(  \sum_{k\in N_{1}}\rho_{jk}+\sum_{k\in N_{2}}%
\rho_{jk}\right)  \right]  \\
= &  \sum_{j\in N_{2}}\sum_{k\in N_{2}}\rho_{kj}-\sum_{j\in N_{2}}\sum_{k\in
N_{2}}\rho_{jk}-\sum_{j\in N_{2}}\sum_{k\in N_{1}}\rho_{jk}\\
= &  -\sum_{j\in N_{2}}\sum_{k\in N_{1}}\rho_{jk}\leq0
\end{split}
\end{equation}
Next, summing the right-hand-sides of the constraints
(\ref{follower_conservation}) over $j\in N_{2}$ we get $\sum_{j\in N_{2}}%
N_{j}=N>0$. This contradicts the constraint (\ref{follower_conservation})
leading to the conclusion that the graph $\mathscr{G(\mathbf{s})}$ is
connected. The additional constraints (\ref{self_loop})-(\ref{flow_capacity})
are necessary to ensure that the number of feasible flow vectors
$\mathbf{\rho}$ is finite. Clearly, (\ref{self_loop}) prohibits self-loops
while (\ref{flow_capacity}) prevents an infinite number of solutions where
edges $(i,j)$ in $\mathscr{G(\textbf{s})}$ may take any unbounded flow value
$\rho_{ij}>0$. $\blacksquare$

Observe that $\rho_{ij}>0$ indicates a connection between agents $i$ and $j$.
This can be combined with the constraint $d_{ij}(t)\leq C_{ij}$ to write
$\rho_{ij}(d_{ij}(t)-C_{ij})\leq0$ for all edges $(i,j)$ in
$\mathscr{G(\textbf{s})}$. Moreover, the convex set $F$ can be expressed
through linear constraints. Thus, the optimal formation problem with
connectivity preservation at any fixed $t\in\lbrack0,T]$ becomes a Mixed
Integer Nonlinear Problem (MINLP):
\begin{equation}%
\begin{split}
&  \min_{\mathbf{s}(t),\mathbf{\rho}}\text{ }-H(\mathbf{s}(t),\mathbf{\rho})\\
\text{s.t.}\quad s_{i}(t) &  \in F\subseteq\Omega,\text{ \ }i=0,\ldots,N\\
\sum_{i\in\mathscr{N}}\rho_{ij} &  -\sum_{i\in\mathscr{N}_{F}}\rho
_{ji}=1,\text{ \ }j\in\mathscr{N}_{F}\\
\rho_{ij}(d_{ij}(t)-C_{ij}) &  \leq0,\text{ }i\in\mathscr{N},\text{ \ }%
j\in\mathscr{N}_{F}\\
\rho_{ii} &  =0,\text{ \ }i\in\mathscr{N}_{F}\\
\rho_{ij} &  \leq N,\text{ \ }i\in\mathscr{N},\text{ \ }j\in\mathscr{N}_{F}%
\end{split}
\label{objwconstraint}%
\end{equation}
Note that any agent position vector $\mathbf{s}(t)$ specifies a graph at time
$t$. The role of $\mathbf{\rho}$ is in ensuring that this graph is connected
by satisfying the constraints in (\ref{objwconstraint}), thus creating an
optimal formation. However, there is no advance information regarding what the
optimal formation looks like and how the optimal formation changes over time
as the leader moves in a time interval $[0,T]$ unless $H(\mathbf{s}(t))$ is
given some specific structure.

For the remainder of this paper, we will consider the class of coverage
control problems
\cite{SM2011,CM2004,cgc2005,caicedo2008coverage,caicedo2008performing,breitenmoser2010voronoi,Minyi2011,Gusrialdi2011}
which impose a particular structure on $H(\mathbf{s}(t))$. Agents are assumed
to be equipped with some sensing and some communication capabilities. In
particular, we assume that agent $i$'s sensing is limited to a set $\Omega
_{i}(t)\subset\Omega$. For simplicity, we let $\Omega_{i}(t)$ be a circle
centered at $s_{i}(t)$ with radius $\delta_{i}$. Thus, $\Omega_{i}%
(t)=\{x:d_{i}(x,t)\leq\delta_{i}\}$ where $d_{i}(x,t)=\Vert x-s_{i}(t)\Vert$,
the standard Euclidean norm. To further maintain simplicity without affecting
the generality of the analysis, we set $\delta_{i}=\delta$ for all agents. We
define $p_{i}(x,s_{i}(t))$ to be the probability that $i$ detects an event
occurring at point $x$. This function is defined to have the following
properties: $(i)$ $p_{i}(x,s_{i}(t))=0$ if $x\notin\Omega_{i}(t)$, and $(ii)$
$p_{i}(x,s_{i}(t))\geq0$ is a monotonically nonincreasing function of
$d_{i}(x,t)$. The overall \emph{sensing detection probability} is denoted by
$\hat{p}_{i}(x,s_{i}(t))$ and defined as
\begin{equation}
{\label{SensingModel}}\hat{p}_{i}(x,s_{i}(t))=%
\begin{cases}
p_{i}(x,s_{i}(t)) & \text{if}\quad x\in\Omega_{i}(t)\\
0 & \text{if}\quad x\notin\Omega_{i}(t)
\end{cases}
\end{equation}
Note that $\hat{p}_{i}(x,s_{i}(t))$ may not be continuous in $s_{i}(t)$. The
joint detection probability, denoted by $P(x,\mathbf{s}(t))$, captures the
sensing ability of the entire agent team. That is, an event at $x\in\Omega$ is
detected by at least one of the $N$ cooperating agents with probability
$P(x,\mathbf{s}(t))$ is given by%
\begin{equation}
P(x,\mathbf{s}(t))=1-\prod_{i=0}^{N}[1-\hat{p}_{i}(x,s_{i}(t))]\label{jointP}%
\end{equation}
where we assume that agents sense independently of each other. In addition to
sensing, the communication capabilities of agents are defined by their
relative distance: agents $i$ and $j$ can establish a communication link if
$\Vert s_{i}(t)-s_{j}(t)\Vert\leq C$. Thus, in this class of problems a
formation is required to maintain full communication among agents. Finally,
one of the agents, indexed by $0$, is designated as the leader whose position
$s_{0}(t)$ is given.

\begin{figure}[h]%
\begin{tabular}
[c]{cc}%
\begin{minipage}[t]{1.5in} \includegraphics[
	height=1.2in,
	width=1.5in]
	{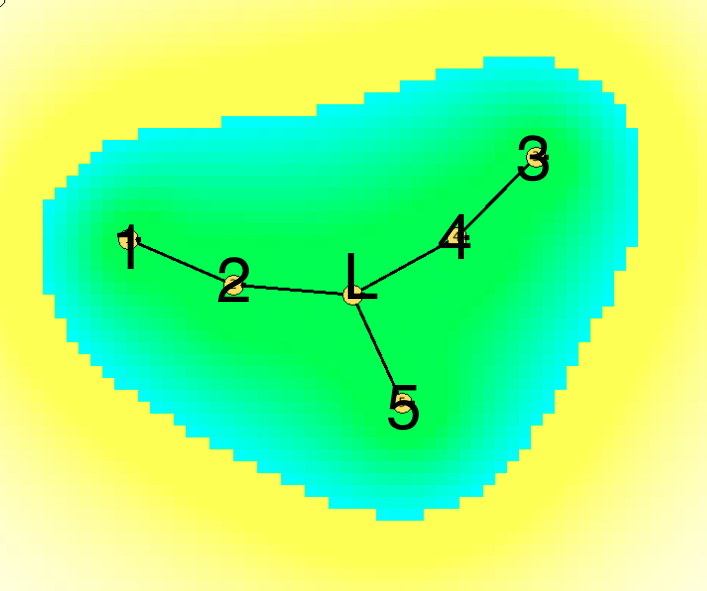} \caption{Optimal formation for 5 followers (numbers) and one leader (L) in a bounded mission space.}\label{fig:ExampleMINLP} \end{minipage} \begin{minipage}[t]{1.5in} \includegraphics[
	height=1.2in,
	width=1.5in]
	{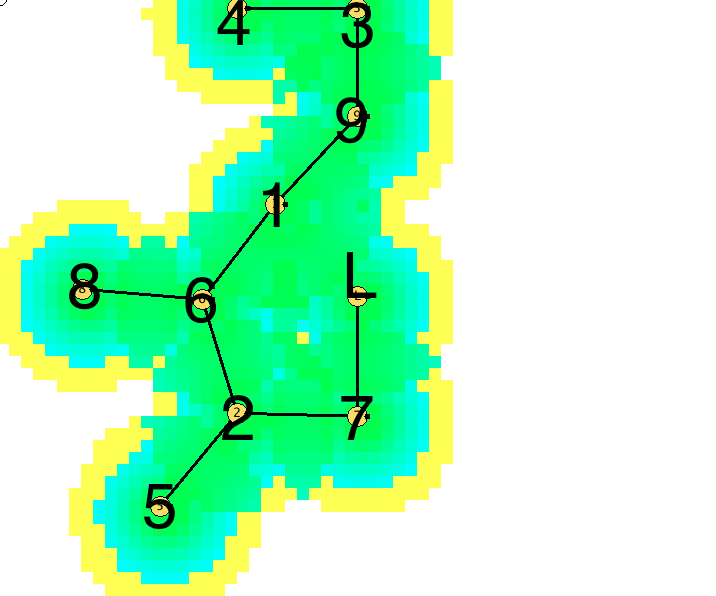} \caption{Optimal formation for 11 followers and a leader. Followers are constrained to the left side of the leader.} \label{fig:MINLPSolution} \end{minipage} &
\end{tabular}
\end{figure}

The objective function for optimal coverage is the same as in \cite{Minyi2011}
except for the presence of a leader whose position is predefined. For any
$x\in\Omega$, the function $R(x):\Omega\rightarrow\mathbb{R}$ captures an a
priori estimate of the frequency of event occurrences at $x$ and is referred
to as an \textquotedblleft event density\textquotedblright\ satisfying
$R(x)\geq0$ for all $x\in\Omega$ and $\int_{\Omega}R(x)dx<\infty$. In this
problem, we assume that the event density is a constant for any $x\in\Omega$.
We are interested in maximizing the total detection probability over the
mission space $\Omega$:
\begin{equation}
\max_{\mathbf{s}(t)}\text{ }H(\mathbf{s}(t))=\int_{\Omega}R(x)P(x,\mathbf{s}%
(t))dx \label{covcontrolH}%
\end{equation}
so that the objective in (\ref{objwconstraint}) is $H(\mathbf{s}%
(t),\mathbf{\rho})=\int_{\Omega}R(x)P(x,\mathbf{s}(t))dx$. Figures
\ref{fig:ExampleMINLP} and \ref{fig:MINLPSolution} show optimal formation
examples obtained by solving (\ref{objwconstraint}) at time $t$ with
$s_{0}(t)$ located at the center of the mission space.

A solution of this MINLP is computationally costly so that it is not realistic
to expect re-solving it over the course of a mission $t\in\lbrack0,T]$ as the
leader moves. However, it is not always necessary to repeatedly solve this
problem over $[0,T]$. Theorem 2 presents a condition under which we only need
to solve the problem at $t=0$. This simply formalizes the rather obvious fact
that if no new constraints (e.g., obstacles) are encountered over $t\in(0,T]$,
then the optimal formation at $t=0$ can be preserved by maintaining fixed
relative positions for all agents.

\textbf{Theorem 2} Let $\mathbf{s}(0)$ be an optimal solution of problem
(\ref{objwconstraint}) at $t=0$ and assume that $\Omega_i(t) \subset F, i \in \mathscr{N}$ and that $s_{0}(t)$ is known to all
followers for all $t\in(0,T]$. If $s_{i}(t)=s_{i}(0)+s_{0}(t)-s_{0}(0)$,
$i\in\mathscr{N}_{F}$, then $\mathbf{s}(t)$ maximizes $H(\mathbf{s}(t))$ in
(\ref{covcontrolH}).

\textbf{Proof:} Let us introduce a local polar coordinate system for each
agent $i$, so that the origin of $i$'s local coordinate system is $s_{i}$ and
the axes are parallel to those in the mission Cartesian coordinate system.
Given any point $x=(x_{x},x_{y})\in F$, let $l=(r_{i},\theta_{i})$ be the
polar coordinates in $i$'s local coordinate system. Then, the transformation
that maps $(r_{i},\theta_{i})$ onto the global coordinate system is
$x=s_{i}(t)+[r_{i}\cos\theta_{i}$ $r_{i}\sin\theta_{i}]^{T}$. Upon switching
to this local coordinate system, the sensing probability becomes
$p_{i}(x,s_{i}(t))=p_{i}(r_{i})$ if $r_{i}<\delta$. Since $\Omega_{i}(t)\in F$
for all $t\in\lbrack0,T]$, the local sensing range of $s_{i}(t)$, which is
denoted by $\Omega_{i}^{L}=\{(r_{i},\theta_{i}):r_{i}\leq C,0\leq\theta
_{i}\leq2\pi\}$, is time-invariant. Therefore, recalling (\ref{jointP}), the
objective function in (\ref{covcontrolH}) is
\begin{equation}%
\begin{split}
H(\mathbf{s}(t))  &  =\int_{\Omega}R(x)P(x,\mathbf{s}(t))dx\\
&  =\int_{\bigcup_{i=0}^{N}\Omega_{i}(t)}R(x)P(x,\mathbf{s}(t))dx\\
&  =\int_{\bigcup_{i=0}^{N}\Omega_{i}(t)}R(x)\{1-\prod_{i=0}^{N}%
[1-p_{i}(x,s_{i}(t))]\}dx\\
&  =\int_{\bigcup_{i=0}^{N}\Omega_{i}^{L}}r_{i}R(r_{i},\theta_{i}%
)\{1-\prod_{i=0}^{N}[1-p_{i}(r_{i})]\}dr_{i}d\theta_{i}%
\end{split}
\label{covHLocalCoordinate}%
\end{equation}
so that the objective function value remains fixed for any $t\in\lbrack0,T]$.
Since for any agents $i$ and $j$, by assumption, $s_{i}(t)-s_{j}%
(t)=s_{i}(0)+s_{0}(t)-s_{0}(0)-\left(  s_{j}(0)+s_{0}(t)-s_{0}(0)\right)
=s_{i}(0)-s_{j}(0)$, and $\mathbf{s}(0)$ is an optimal solution of
(\ref{objwconstraint}), it follows that $\mathscr{G}(\mathbf{s}(0))$ is
connected, therefore, $\mathscr{G}(\mathbf{s}(t))$ is also connected and we
conclude that $\mathbf{s}(t)$ maximizes $H(\mathbf{s}(t))$. $\blacksquare$

The implication of Theorem 2 is that when a mission space has no obstacles in
it or the leader follows a trajectory where no obstacles are encountered by
any agent, our problem is reduced to one of ensuring that all agents
accurately track the leader's trajectory. We may discretize time so that
agents update their locations at $0<t_{1}<\cdots<t_{K}=T$. Assuming that
problem (\ref{objwconstraint}) is solved at $t=0$, an optimal formation is
obtained and we subsequently strive to maintain this formation until a
significant \textquotedblleft event\textquotedblright\ occurs such as an agent
failure, a change in objective function $H(\mathbf{s}(t))$, or encountering
obstacles; at such a point, some amount of reconfiguration is required while
still aiming to maximize $H(\mathbf{s}(t))$.

\section{Optimal Dynamic Formation Control in a Mission Space with Obstacles}

We have thus far solved an optimal dynamic formation problem with connectivity
constraints in a convex feasible space $F$ by solving a MINLP. However, this
method may fail when $F$ is non-convex, e.g., when $F$ cannot be described
through linear or nonlinear constraints. In this section, we address the
optimal dynamic formation problem in a mission space with obstacles, thus
considering a non-convex feasible space.

We model the obstacles as $m$ non-self-intersecting polygons denoted by
$M_{j}$, $j=1,\ldots,m$. The interior of $M_{j}$ is denoted by $\mathring
{M_{j}}$, so that the overall feasible space is $F=\Omega\setminus
(\mathring{M_{1}}\cup\ldots\cup\mathring{M_{m}})$, i.e., the space $\Omega$
excluding all interior points of the obstacles. In this setting, we seek to
ensure the following two requirements. First, the distance between two
connected agents must be $\leq C$. We define $c_{1}(s_{i},s_{j})$ to indicate
whether this requirement is satisfied:
\begin{equation}
c_{1}(s_{i},s_{j})=\left\{
\begin{array}
[c]{cl}%
1 & \Vert s_{i}-s_{j}\Vert\leq C\\
0 & \text{otherwise}%
\end{array}
\right.
\end{equation}
Second, the connected agents are required to have a line of sight with respect
to each other. We define $c_{2}(s_{i},s_{j})$ to indicate this requirement:
\begin{equation}
c_{2}(s_{i},s_{j})=\left\{
\begin{array}
[c]{cl}%
1 & \alpha s_{i}+(1-\alpha)s_{j}\in F\text{ for all }\alpha\in\lbrack0,1]\\
0 & \text{otherwise}%
\end{array}
\right.
\end{equation}
Agents $i$ and $j$ satisfying $c_{1}(s_{i},s_{j})=1$ as well as $c_{2}%
(s_{i},s_{j})=1$ are referred to as \emph{connected}. We also define
$c(s_{i},s_{j})=c_{1}(s_{i},s_{j})c_{2}(s_{i},s_{j})$.

A version of this connectivity preservation problem was addressed in
\cite{Minyi2011}, where agents are required to remain connected with a fixed
base while at the same time maximizing the objective function in
(\ref{covcontrolH}). A gradient-based algorithm, termed \emph{Connectivity
Preservation Algorithm} (CPA), was developed for agent position updating and
it was shown that, given an initially connected network and if only one agent
updates its position at any given time, the CPA preserves connectivity. The
algorithm is applied iteratively over one agent at a time and it converges to
a (generally local) optimum. The CPA exploits the existence of distributed
optimization algorithms for optimal coverage to attain optimal agent locations
while also preserving connectivity to a base (details on the CPA and its
complexity are provided in \cite{Minyi2011}). 

\begin{figure}[ptb]
\centering
\includegraphics[
height=2.1in,
width=3.4in]
{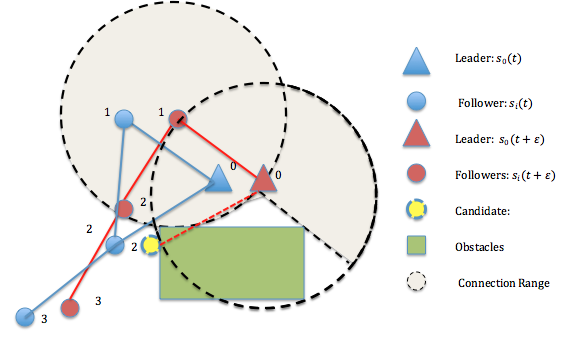}\caption{An example of a connected network at $t$ and
constructed connected network by Algorithm 1 at $t+\epsilon$.}%
\label{fig: Algorithm}%
\end{figure}

Our approach here is to take advantage of the CPA. In our problem, however,
the conditions for applying the CPA do not generally hold; this is because the
leader's motion does not take connectivity with its neighbors into account and
the presence of an obstacle, for example, may cause it to disconnect from one
or more followers. This is illustrated in Fig. \ref{fig: Algorithm}: At time
$t$, the agent network shown (represented by three blue circles and a blue
triangle as the leader) is connected. At $t+\epsilon$, the leader (triangle)
moves to $s_{0}(t+\epsilon)$ and if agent 2 moves to the point shown in yellow
(as expected by Theorem 2), then it becomes disconnected from the leader
because of the obstacle present. We propose an algorithm next to construct a connected graph, which
may no longer be optimal in the sense of problem (\ref{objwconstraint}) but it
does provide a valid initial condition for invoking the CPA described above
(this is illustrated in Fig. \ref{fig: Algorithm} as the solid red graph).
This immediately allows us to iteratively apply the CPA so as to obtain a new
(locally optimal) formation.

Clearly, it is also possible to invoke (\ref{objwconstraint}) as soon as a
formation reconfiguration is needed. However, the set $F$ is no longer convex
and the computational complexity of this problem makes it infeasible for the
on-line adaptation required, whereas the approach we propose and the use of
the CPA render this process computationally manageable. In particular, whereas
the MINLP is generally NP hard, in the CPA each agent $i$ determines its new
position through a gradient-based scheme using only its neighbor set and its downstream and upstream agent sets relative to the leader (formally defined in the next section). When the
number of agents increases, note that the the number of neighbors of $i$ may not be affected. The overall increase in complexity is linear in the network size.

\begin{figure}[ptb]
\centering
\includegraphics[
height=2in,
width=2.2in]
{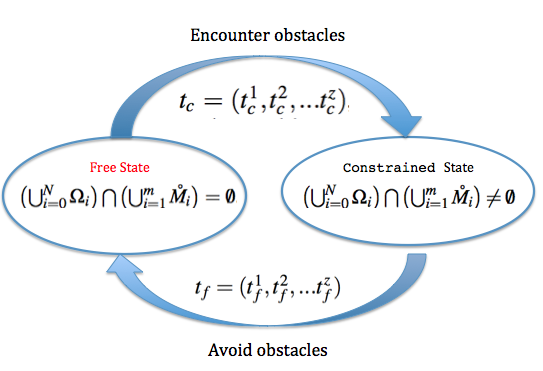} .\caption{Two states of the agents network and the
transition time points between the two states.}%
\label{fig:stateTransition}%
\end{figure}

Before proceeding, we identify the precise instants when formation
reconfiguration is necessary due to obstacles encountered by agents as the
mission unfolds over $[0,T]$. We define two states that the agent team can be
in: $(i)$ The \emph{constrained} state occurs when the sensing capability of
an agent is hindered by an obstacle, captured by the condition $\left(
\bigcup_{i=0}^{N}\Omega_{i}\right)  \bigcap\left(  \bigcup_{i=1}^{m}%
\mathring{M_{i}}\right)  \neq\emptyset$, i.e., the intersection of the sensed
part of $\Omega$ and the set of interior points of any obstacle is not empty,
and $(ii)$ The \emph{free} state corresponding to $\left(  \bigcup_{i=0}%
^{N}\Omega_{i}\right)  \bigcap\left(  \bigcup_{i=1}^{m}\mathring{M_{i}%
}\right)  =\emptyset$. Thus, the interval $[0,T]$ is partitioned into free and
constrained intervals with transitions at times $t_{f}^{0}<t_{c}^{1}<t_{f}%
^{1}<...<t_{c}^{i}<t_{f}^{i}<...t_{f}^{z}<T$. This is described in Fig.
\ref{fig:stateTransition}. Next, we consider how to generate optimal
formations over different alternating intervals $[t_{f}^{k},t_{c}^{k+1})$ and
$[t_{c}^{k+1},t_{f}^{k+1})$. 

\subsection{Optimal formation control in free states}

When the agent network enters a free state at time $t_{f}^{k},$ $k=0,\ldots
,z$, since $\left(  \bigcup_{i=0}^{N}\Omega_{i}(t)\right)  \bigcap\left(
\bigcup_{i=1}^{m}\mathring{M_{i}}\right)  =\emptyset$ for all $t \in [t_f^k,t_c^{k+1})$ and $F=\Omega
\setminus(\mathring{M_{1}}\cup\ldots\cup\mathring{M_{m}})$, so $\Omega_{i}(t)\in F$ for any $i$ over $t \in [t_f^k,t_c^{k+1})$, the optimal formation is maintained based on Theorem 2.

\subsection{Optimal formation control in constrained states}

We begin this subsection with some additional notation and definitions. Given
a connected graph $\mathscr{G}(\mathbf{s})$, we have defined a loop-free path
connecting agent $i$ to the leader as $\pi_{i}=\{0,\ldots,a,b,\ldots,i\}$, an
ordered set where neighboring agents are connected; we have also defined
$\Pi_{i}$ to be the set of all possible paths connecting $i$ to the leader.
Let $\pi_{i,k}$ be the $k$th path in $\Pi_{i}$ and we use $\pi_{i,k}^{j}$ to
denote the $j$th element in $\pi_{i,k}$. Let $\mathscr{D}_{i}=\cup_{j,k}%
w_{i}(\pi_{j,k})$ be the set of agents \emph{downstream}  from $i$ (further
away from the leader 0) where
\begin{equation}
w_{i}(\pi_{j,k})=\left\{
\begin{array}
[c]{cl}%
\pi_{j,k}^{l+1} & \text{ if }i\in\pi_{j,k},\text{ }i\neq j\text{ and }%
i=\pi_{j,k}^{l}\\
\emptyset & \text{otherwise}%
\end{array}
\right.
\end{equation}
We also define the set of \emph{upstream} agents from $i$ as $\mathscr{U}_{i}%
=\{j:i\in\mathscr{D}_{j},j\in{0,\dots,N}\}$.

The length of a path $\pi_{i,k}$ is defined as $\Psi(\pi_{i,k})=\sum
_{l=1}^{|\pi_{i,k}|-1}\Vert s_{\pi_{i,k}^{l}}-s_{\pi_{i,k}^{l+1}}\Vert$, where
$|\pi_{i,k}|$ is the cardinality of $\pi_{i,k}$. For agent $i$, the shortest
path connected to the leader is
\[
\pi_{i}^{\ast}=\arg\min_{\pi_{i,k}\in\Pi_{i}}\Psi(\pi_{i,k})
\]
For example in Fig. \ref{fig: Algorithm}, in the path $\pi_{3,1}=\{0,2,3\}$,
we have 3 $\in\mathscr{D}_{2}$, 0 $\in\mathscr{U}_{2}$, $\Psi(\pi_{3,1})=\Vert
s_{0}-s_{2}\Vert+\Vert s_{2}-s_{3}\Vert$; for the path $\pi_{3,2}%
=\{0,1,2,3\}$, we have $\Psi(\pi_{3,2})=\Vert s_{0}-s_{1}\Vert+\Vert
s_{1}-s_{2}\Vert+\Vert s_{2}-s_{3}\Vert$. Therefore, $\pi_{3}^{\ast}=\pi
_{3,1}$ is the shortest path from agent $3$ to the leader.

Let $\pi_{i}$ and $\pi_{j}$ be two paths. Then, we define $\pi_{i}+\pi
_{j}=\{\pi_{i},\pi_{k}\}$, where $\pi_{k}=\pi_{j}\setminus\pi_{i}$, as an
ordered set. Note that $\pi_{i}+\pi_{j}$ is generally different from $\pi
_{j}+\pi_{i}$ because of the order involved. Given a connected graph
$\mathscr{G}(\mathbf{s})$, We define
\begin{equation}
Q(\mathscr{G}(\mathbf{s}))=\pi_{1}^{\ast}+\ldots+\pi_{N}^{\ast}\label{Qdef}%
\end{equation}
to be an ordered set containing a permutation of the agent set $\{0,1,...,N\}$
constructed so as to start with the shortest path $\pi_{1}^{\ast}$ from $0$ to
agent $1$, followed by $\pi_{2}^{\ast}\setminus\pi_{1}^{\ast}$ and so on. It
immediately follows from this construction that the first element of
$Q(\mathscr{G}(\mathbf{s}))$ is $0$ and that $|Q(\mathscr{G}(\mathbf{s}%
))|=N+1$. Therefore, we can rewrite $Q(\mathscr{G}(\mathbf{s}))$ as
\[
Q(\mathscr{G}(\mathbf{s}))=\{0,q_{2},\ldots,q_{N+1}\}
\]
where $q_{j}\in\mathscr{N}_{F},j=2,\ldots,N+1$. For example, in Fig.
\ref{fig: Algorithm}, at time $t$, $Q(\mathscr{G}(\mathbf{s}(t)))=\{0,1,2,3\}$%
. We show next that $Q(\mathscr{G}(\mathbf{s}))$ has the following property
regarding the order of its elements.

\textbf{Lemma 1 } If $q_{i}$ is the $i$th element of $Q(\mathscr{G}(\mathbf{s}%
))$ constructed from a connected graph $\mathscr{G}(\mathbf{s})$, then there
exists $q_{j}\in\mathscr{U}_{q_{i}}$ such that $q_{j}$ is the $j$th element of
$Q(\mathscr{G}(\mathbf{s}))$, and $j<i$ for all $q_{i}\in\mathscr{N}_{F}$.

\textbf{Proof:} If for all $q_{j}\in\mathscr{U}_{q_{i}}$, $j>i$, we cannot
find a subset of $Q(\mathscr{G}(\mathbf{s}))$ that includes $\{q_{j},q_{i}\}$,
$q_{j}\in\mathscr{U}_{q_{i}}$, then there is no path connected to $q_{i}$.
This contradicts the assumption that $Q(\mathscr{G}(\mathbf{s}))$ is
constructed from a connected graph. $\blacksquare$

We also define a projection of $x\in\mathbb{R}^{2}$ on a set $A\in
\mathbb{R}^{2}$ as%
\[
P_{A}(x)=\arg\min_{y\in A}\Vert x-y\Vert
\]
Next, let $\mathscr{Y}(s_{i})=\{y:y\in\mathbb{R}^{2},$ $c(s_{i},y)=1)$.
Recalling the definition of $c(\cdot,\cdot)$, $\mathscr{Y}(s_{i})$ is the set
of points with which $s_{i}$ can establish a connection. For any subset of
agents $\mathscr{V}\subset\mathscr{N}$, let $\Sigma(\mathscr{V})=\bigcup
_{i\in\mathscr{V}}\mathscr{Y}(s_{i})$ be the union of all connection regions
for agents in $\mathscr{V}$. For example, in Fig. \ref{fig: Algorithm}, the
grey area is $\Sigma(\mathscr{V})$ for $\mathscr{V}=\{0,1\}$ at time
$t+\epsilon$.

We are now ready to deal with the situation where the formation is in a
constrained state and may lose connectivity at time $t+\epsilon$ given that
the graph $\mathscr{G}(\mathbf{s}(t))$ is connected. In particular, suppose
that when the leader is about to move to $s_{0}(t+\epsilon)$ and informs the
followers, at least one of the agents will lose connectivity with the
formation. Our task is to obtain an optimal formation at $t+\epsilon$ and this
is accomplished in two steps: $(i)$ Construct a connected graph
$\mathscr{G}(\mathbf{s}(t+\epsilon))$ for time $t+\epsilon$, and $(ii)$ Use
this connected graph $\mathscr{G}(\mathbf{s}(t+\epsilon))$ as an input to
invoke the CPA. Step $(i)$ is crucial because of the fact that the CPA relies
on an initially connected graph before it can be executed to seek (locally)
optimal agent locations which still preserve connectivity. This first step is
carried out by constructing a connected graph through Algorithm
\ref{alg: ConGraph}. \begin{algorithm}
\caption{ Connected Graph Construction Algorithm}
\label{alg: ConGraph}
\textbf{Input}: Graph $\mathscr{G}(\mathbf{s}(t))$, $s_0(t+\epsilon)$ \\
\textbf{Output}: Graph $\mathscr{G}(\mathbf{s}(t+\epsilon))$ \\
\textbf{Initialization:} $\mathscr{U}_i, \mathscr{D}_i$ for $i \in \mathscr{N}$, $\mathscr{V}=\{ 0 \}$,    $Q(\mathscr{G}(\mathbf{s}(t)))=\{0, q_2,\ldots, q_{N+1} \}$ using (\ref{Qdef}) \\
\textbf{For} agent $i=q_j, j=2,\ldots, N+1$ \\
Do the following procedure:
\begin{algorithmic}[1]
\STATE{\label{Candidate}} Generate a candidate next location for $i$: $\hat{s}_i=s_i(t)+\Delta_L$.
\STATE{\label{UpstreamCheck}} For all agents $v \in \mathscr{U}_i \bigcap \mathscr{V}$, if $c(\hat{s}_i, s_v(t+\epsilon)) = 0$, go to Step \ref{UpstreamProj}; else, go to Step \ref{Update}.
\STATE{\label{UpstreamProj}} Project $s_i$ onto $\Sigma(\mathscr{U}_i \bigcap \mathscr{V})$. Set $\hat{s_i} = P_{\Sigma(\mathscr{U}_i \bigcap \mathscr{V})}(s_i)$.
\STATE{\label{Update}}  Set $s_i(t+\epsilon) = \hat{s}_i$.
\STATE{\label{UpdateI}}  Add $i$ to $\mathscr{V}$
\end{algorithmic}
\textbf{End}
\end{algorithm}We use $\Delta_{L}(t)=s_{0}(t+\epsilon)-s_{0}(t)$ to denote the
position change vector of the leader from $t$ to $t+\epsilon$, where we assume
that followers have the $\Delta_{L}(t)$ information available at $t$.

\textbf{Theorem 3 } $\mathscr{G}(\mathbf{s}(t+\epsilon))$ obtained by
Algorithm \ref{alg: ConGraph} is connected.

\textbf{Proof}: Since $\mathscr{G}(\mathbf{s}(t))$ is connected,
$\mathscr{U}_{i}\neq\emptyset$ for $i\in\mathscr{N}_{F}$. We then use
induction to prove that the graph constructed by agents in $\mathscr{V}$
remains connected at Step \ref{UpdateI} in every iteration. Initially,
$\mathscr{V}=\{0\}$ which is connected. Next, assuming there are $n$ agents in
$\mathscr{V}$ and the graph they form is connected, we will prove that after
adding the $(n+1)$th agent, say $i$, the graph remains connected.

The addition of $i$ to $\mathscr{V}$ occurs at Step 5. There are two possible
sequences for reaching this step: \ref{Candidate}-\ref{UpstreamCheck}%
-\ref{Update} and \ref{Candidate}-\ref{UpstreamCheck}-\ref{UpstreamProj}%
-\ref{Update}. At Step \ref{UpstreamCheck}, $\mathscr{U}_{i}\bigcap
\mathscr{V}\neq\emptyset$ because of the property of $Q(\mathscr{G}(\mathbf{s}%
))$ in Lemma 1. It follows that before $i$ performs the procedure, there is at
least one upstream agent in $\mathscr{V}$. In the \ref{Candidate}%
-\ref{UpstreamCheck}-\ref{Update} sequence, there exists some $m\in
\mathscr{V}\cap\mathscr{U}_{i}$ such that $c(\hat{s}_{i},s_{m}(t+\epsilon
))=1$. Therefore, all agents in $\mathscr{V}$ including $i$ will be connected.
In the \ref{Candidate}-\ref{UpstreamCheck}-\ref{UpstreamProj}-\ref{Update}
sequence, at Step \ref{UpstreamProj}, agent $i$'s position is projected onto
the connection ranges of all $v\in\mathscr{V}\cap\mathscr{U}_{i}$. It follows
that the graph formed by agents in $\{\mathscr{V},i\}$ is connected. Step 5
adds agents to $\mathscr{V}$ one by one until $\mathscr{V}=\mathscr{N}$,
therefore, the graph $\mathscr{G}(\mathbf{s}(t+\epsilon))$ is connected.
$\blacksquare$

Obviously, Algorithm \ref{alg: ConGraph} does not provide a unique way to
construct a connected graph. For example, the formation could be adjusted to a
line or a star configuration with $s_{0}(t+\epsilon)$ as the center of the
star. However, this would entail a major formation restructuring whereas in
Algorithm \ref{alg: ConGraph} we seek to retain the \emph{closest possible
formation} to the original (optimal) one by setting candidate locations as
seen in Step \ref{Candidate}. If such a candidate is not feasible, then the
agent will move a minimal distance (in the projection sense) to be connected.

Once step $(i)$ above is completed by obtaining this connected graph
$\mathscr{G}(\mathbf{s}(t))$, step $(ii)$ is performed by invoking the CPA to
optimize the agent locations within the new formation. Clearly, once obstacles
are cleared and the agent team re-enters a free state (see Fig.
\ref{fig:stateTransition}), we may revert to the original optimal formation.

\section{Simulation Results}

In this section, we provide a simulation example illustrating what the optimal
formation maximizing coverage in a mission space with obstacles looks like and
how it changes at some significant instants.

We choose the event density functions to be uniform, i.e., $R(x)=1$. The
mission space is a $60\times50$ rectangle. The distance constraint is $C=10$
and the sensing range of each agent is $\delta=8$. At every step, the leader
moves to the right one distance unit per unit of time. The mission space is
colored from dark to lighter as the joint detection probability decreases (the
joint detection probability is $\geq0.50$ for green areas, and near zero for
white areas). The leader (labeled \textquotedblleft L\textquotedblright) moves
along a predefined trajectory (the purple dashed line). There are 8 followers,
indicated by numbers, which are restricted to locations on the left
side of the leader during any movement.

\begin{figure}[h]%
\begin{tabular}
[c]{cc}%
\begin{minipage}[t]{1.6in} \includegraphics[
	height=1.2in,
	width=1.6in]
	{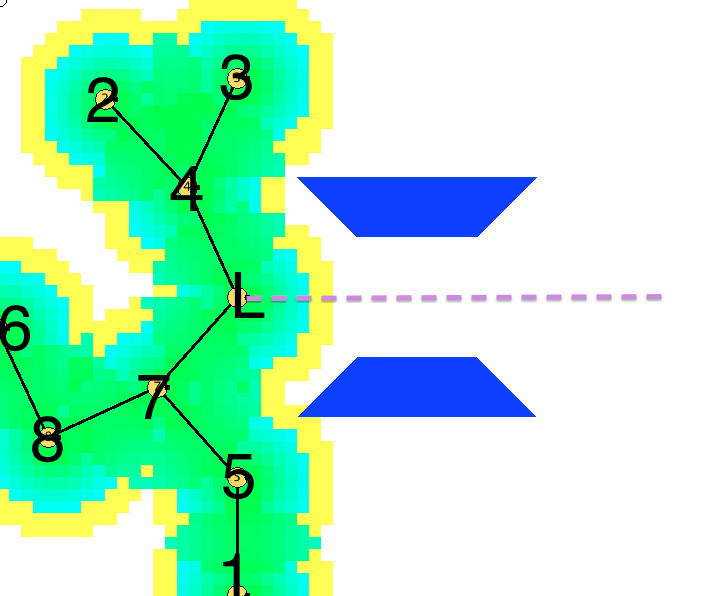} \caption{At $t=0$, optimal formation from MINLP, $H(\mathbf{s})=741.5$} \label{fig:MINLPStart} \end{minipage} \begin{minipage}[t]{1.6in} \includegraphics[
	height=1.2in,
	width=1.6in]
	{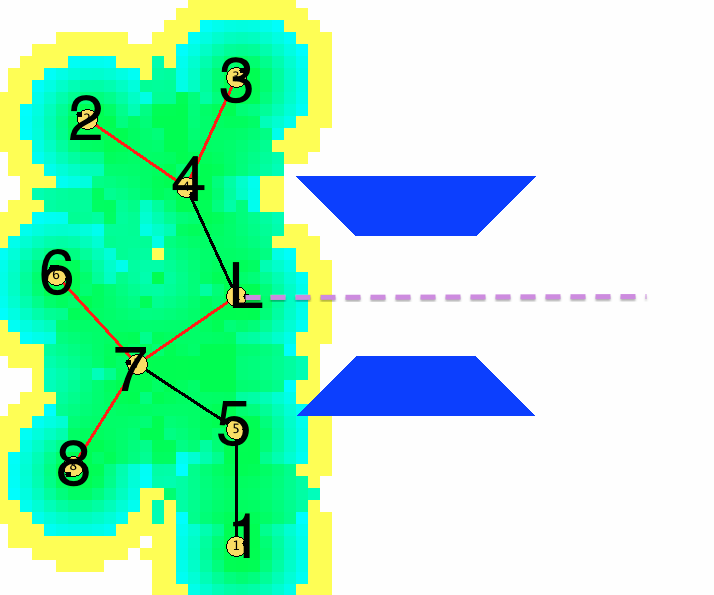} \caption{At $t=0$, optimal formation improved by CPA of Fig. \ref{fig:MINLPStart} with $H(\mathbf{s})=816.7$} \label{fig:CPAFromMINLP} \end{minipage} &
\\
\begin{minipage}[t]{1.6in} \includegraphics[
	height=1.2in,
	width=1.6in]
	{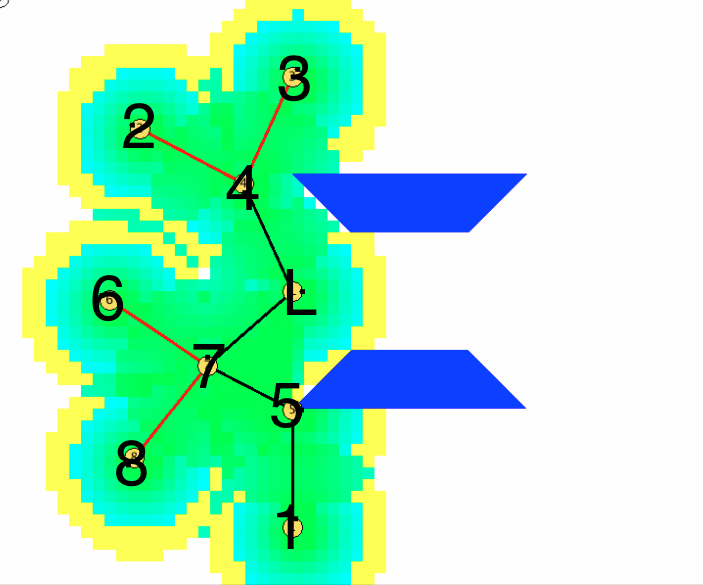} \caption{At $t=5$, agent 5 needs projection in Step 3 of Algorithm 1} \label{fig:firstProjection} \end{minipage} \begin{minipage}[t]{1.6in} \includegraphics[
	height=1.2in,
	width=1.6in]
	{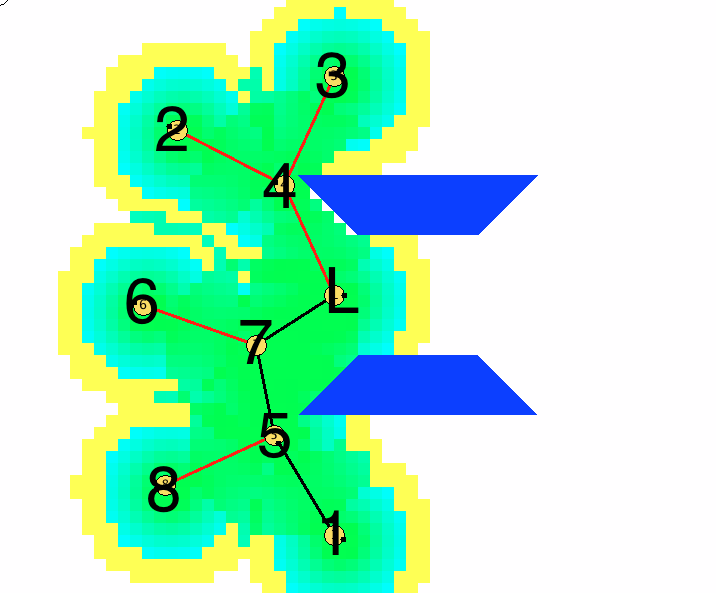} \caption{At $t=7$, apply CPA after projection in \ref{fig:firstProjection}, the structure of the tree doesn't change} \label{fig: CPAProjection} \end{minipage} &
\end{tabular}
\end{figure}

\begin{figure}[h]%
\begin{tabular}
[c]{cc}%
\begin{minipage}[t]{1.6in} \includegraphics[
	height=1.2in,
	width=1.6in]
	{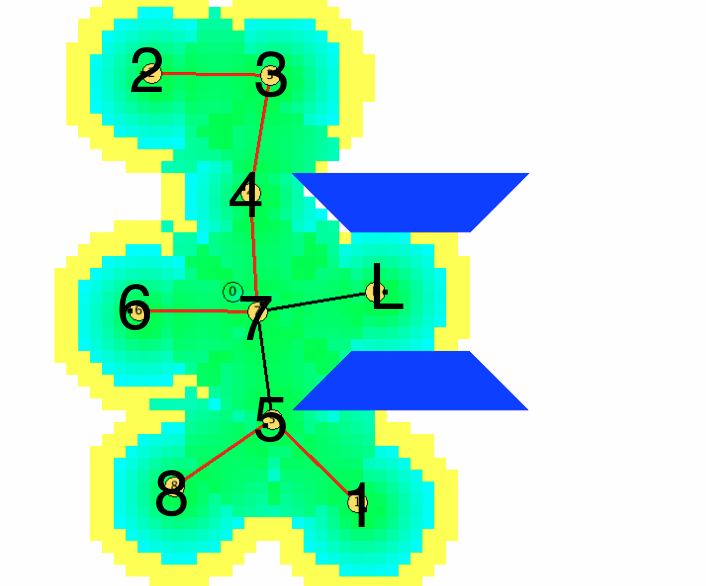} \caption{At $t=12$, the structure of the tree changes} \label{fig:Reconfiguration} \end{minipage} \begin{minipage}[t]{1.6in} \includegraphics[
	height=1.2in,
	width=1.6in]
	{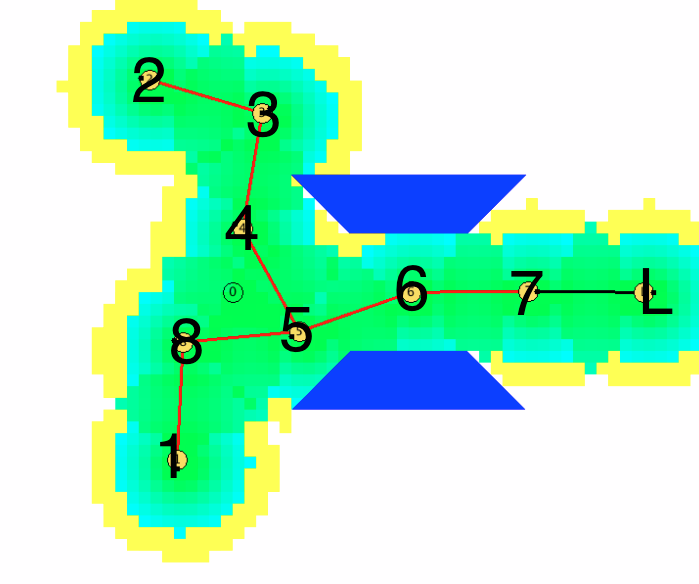} \caption{At $t=35$, the end of the mission} \label{fig:MissionEnd} \end{minipage} &
\\
\begin{minipage}[t]{1.6in} \includegraphics[
	height=1.4in,
	width=1.6in]
	{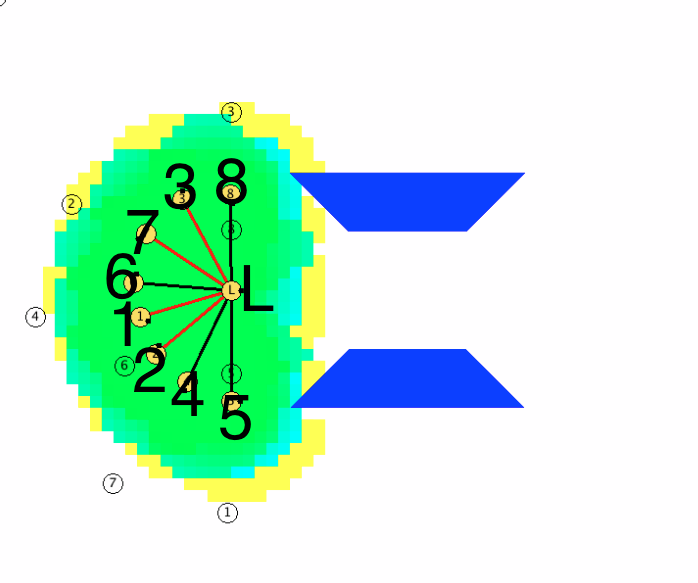} \caption{A star-like connected graph} \label{fig:starStart} \end{minipage} \begin{minipage}[t]{1.6in} \includegraphics[
	height=1.3in,
	width=1.6in]
	{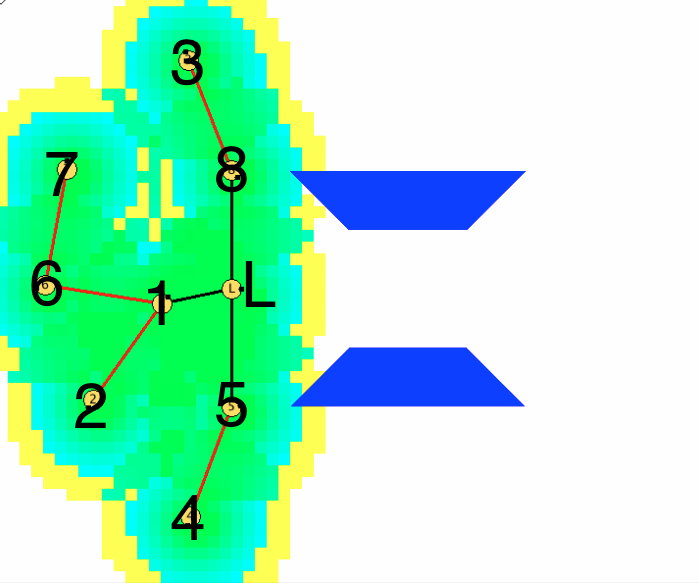} \caption{Apply CPA from \ref{fig:starStart}, $H(\mathbf{s})=781.1$} \label{fig:CPAStart} \end{minipage} &
\end{tabular}
\end{figure}

Figures \ref{fig:MINLPStart}-\ref{fig:MissionEnd} show snapshots of the
process at selected events of interest over $[0,T]$. Figure
\ref{fig:MINLPStart} shows the initial configuration at $t=0$, where the agent
team is located in a convex feasible space. As shown in Sec. III, in this
case, the optimal formation can be obtained by solving a MINLP
\cite{Bussieck2014}. In the results shown, we have used TOMLAB, a MATLAB-based
optimization solver. For the non-convex objective function defined in
(\ref{covcontrolH}), the solution is usually a local maximum; we sought to
find the best local (possibly global) optimum possible by implementing a
multi-start algorithm on the solver. This is done at the start of the mission,
when an off-line computationally intensive procedure is possible. Moreover,
this local maximum can be improved by applying the CPA; in fact, in this
example the use of the CPA led to an improvement from $H(\mathbf{s})=741.5$ to
$H(\mathbf{s})=816.7$, as shown in Fig. \ref{fig:CPAFromMINLP}. Thus, in
general, supplying the CPA with an initial connected graph obtained by solving
the MINLP enables it to converge to a better value. For example, Fig.
\ref{fig:CPAStart} is a local maximum attained by starting with a star-like
connected graph shown in Fig. \ref{fig:starStart} with the objective function
value $H(\mathbf{s})=781.1$ (although this is still worse than the value in
Fig. \ref{fig:CPAFromMINLP}).

In the time interval $[0,5]$, the formation is maintained. At $t=5$, agent 5
is located at a vertex of an obstacle and will therefore lose connectivity as
the leader moves to the next step at $t=6$. At this point, agent 5 will
determine its next position $s_{5}(6)$ by applying a projection at Step
\ref{UpstreamProj} of Algorithm \ref{alg: ConGraph}. Note that only agent 5
needs to perform this projection, rather than the whole team of agents, hence
the computational effort is minimal. Figure \ref{fig: CPAProjection} captures
the optimal formation following Fig. \ref{fig:firstProjection}.

Observe that over the period $[0,12)$, although the optimal formation remains
a tree, it is no longer the same as the original one. However, for each agent
$i$, its downstream node set $\mathscr{D}_{i}$ and upstream node set
$\mathscr{U}_{i}$ remain unchanged. At $t=12$, clearly, the structure of the
formation has been changed. This is a consequence of either the projection
step in Algorithm 1 or the CPA. At the end of the mission at $t=35$, the
formation is shown in Fig. \ref{fig:MissionEnd}. The agents seek to form a
line to go through the narrow region of the mission space while at the same
time maximizing coverage. During the remaining interval $[12,35]$, the process
is similar to what is seen over $[5,12]$.

As we pointed out in the last section, constructing a connected graph can be
accomplished in a variety of ways. As shown in Fig. \ref{fig:starStart}, a
star-like graph is an inferior formation to that of Fig.
\ref{fig:CPAFromMINLP}; this is expected since the latter was obtained
specifically to maximize the objective function in (\ref{covcontrolH}). In
addition, a reconfiguration process as shown in Fig. \ref{fig:CPAStart}
requires agents to move longer distances, hence consuming more energy.

\section{Conclusions and future work}

We have addressed the issue of optimal dynamic formation of multi-agent
systems in mission spaces with obstacles. When the agent team is in a free
state (no obstacles in the mission space affecting them), a locally optimal
solution of a MINLP can provide an initial formation that agents maintain or
it is a good initial point for using the CPA (developed in prior work
\cite{Minyi2011}) to obtain a better local optimum. When the feasible space is
non-convex and connectivity is lost, we have developed an algorithm to
construct a connected graph as an input for the CPA while seeking to maintain
the original formation with minimal effort.

Future work aims at investigating optimal dynamic formation control for more
general classes of objective functions, beyond the coverage control problem.

\bibliographystyle{IEEEtran}
\bibliography{research}

\end{document}